\documentclass[english,aps,prb,amsmath,tightenlines,notitlepage]{revtex4-1}
\usepackage{mathptmx}

\usepackage[T1]{fontenc}
\usepackage[latin9]{inputenc}
\usepackage{float}
\usepackage{amsmath}
\usepackage{amssymb}
\usepackage{graphicx}

\makeatletter

\providecommand{\tabularnewline}{\\}

\makeatother

\usepackage{babel}
\begin{document}
\title{A fast algorithm for computing the Boys function}
\author{Gregory Beylkin}
\email{beylkin@colorado.edu}

\affiliation{Department of Applied Mathematics, University of Colorado, Boulder,
CO 80309, USA}
\author{Sandeep Sharma}
\email{sanshar@gmail.com}

\affiliation{Department of Chemistry, University of Colorado, Boulder, CO 80302,
USA}
\begin{abstract}
We present a new fast algorithm for computing the Boys function using
a nonlinear approximation of the integrand via exponentials. The resulting
algorithms evaluate the Boys function with real and complex valued
arguments and are competitive with previously developed algorithms
for the same purpose.
\end{abstract}
\maketitle

\section{Introduction}

The Boys function \citep{BOYS:1950}
\begin{equation}
F\left(n,z\right)=\int_{0}^{1}e^{-zt^{2}}t^{2n}dt=\frac{1}{2}\int_{0}^{1}e^{-zs}s^{n-1/2}ds,\label{eq:boys function}
\end{equation}
appears in problems of computing Gaussian integrals, and over the
years, there were many algorithms proposed for its evaluation, see
e.g. \citep{MATSUO:1972,WEIDEM:1994,CAR-POL:1998,PEE-KNI:2020,POL-CAR:2002,MAMEDO:2004,MATHAR:2004,WEI-OCH:2015,MA-MA-LA:2016}.
The Boys function is related to a number of special functions, for
example the error function and the incomplete Gamma function, and
(for pure imaginary argument) to the Fresnel integrals.

It is common (see e.g. \citep{MAMEDO:2004,PEE-KNI:2020}) to use recursion
to compute the Boys function for different $n$. The recursion is
obtained via integration by parts, 
\begin{eqnarray}
F\left(n,z\right) & = & -\frac{1}{2z}\int_{0}^{1}\frac{d}{ds}\left(e^{-zs}\right)s^{n-1/2}ds\nonumber \\
 & = & \frac{n-1/2}{z}F\left(n-1,z\right)-\frac{1}{2z}e^{-z},\label{eq:forward recursion}
\end{eqnarray}
and can be run starting with $n=1$ so that we need to have the value
$F\left(0,z\right)$ or starting from a large $n=n_{\max}$ and going
to $n=1$ 
\begin{equation}
F\left(n-1,z\right)=\frac{z}{n-1/2}F\left(n,z\right)+\frac{1}{2\left(n-1/2\right)}e^{-z},\label{eq:backward recursion}
\end{equation}
so that we need to have the value $F\left(n_{\max},z\right)$. In
order to avoid a loss of accuracy, the choice of which recursion to
use depends on the size $z$ and $n_{\max}$. Iterating recursion
\eqref{eq:forward recursion}, the dominant term expressing $F\left(n,z\right)$
via $F\left(0,z\right)$ is $\prod_{j=1}^{n}\left(j-1/2\right)/z^{n}$.
We set 
\[
z^{*}=\left(\prod_{j=1}^{n}\left(j-1/2\right)\right)^{1/n}
\]
and choose \eqref{eq:forward recursion} when $\left|z\right|\ge z^{*}$
and \eqref{eq:backward recursion} if $\left|z\right|<z^{*}$. For
example, if $n_{\max}=18$ then $z^{*}\approx6.75$. We note that
other choices of the parameter $z^{*}$ are possible.

At each step, recursions \eqref{eq:forward recursion} and \eqref{eq:backward recursion}
require only three multiplications and one addition (since the coefficients
can be computed in advance and stored), so it is hard to obtain a
more efficient alternative if one needs to compute these functions
for a range of $n$, $1\le n\le n_{\max}$. In order to initialize
these recursions, we need fast algorithms for computing $F\left(0,z\right)$
and $F\left(n_{\max},z\right)$. Computing $F\left(0,z\right)$ for
real $z$ is straightforward since 
\begin{equation}
F\left(0,z\right)=\int_{0}^{1}e^{-zt^{2}}dt=\frac{\sqrt{\pi}\mbox{Erf}\left(\sqrt{z}\right)}{2\sqrt{z}}.\label{eq:Boys and Error fnc}
\end{equation}
For a real argument an optimized implementation of the error function
$\mbox{Erf}$ is available within programming languages. For a complex
argument we present an algorithm for computing $F\left(0,z\right)$
using nonlinear approximation of the integrand following the approach
in \citep{BEY-MON:2016}. We obtain a rational approximation of $F\left(0,z\right)$
with an additional exponential factor.

We note that, as a function of complex argument, the Boys function
$F\left(0,z\right)$ can be highly oscillatory. In particular, if
$z$ is purely imaginary, then the Boys function is related to the
Fresnel integrals, 
\[
S\left(y\right)=\int_{0}^{y}\sin\left(\frac{\pi}{2}t^{2}\right)dt,\,\,\,\,\,C\left(y\right)=\int_{0}^{y}\cos\left(\frac{\pi}{2}t^{2}\right)dt,
\]
so that
\begin{equation}
C\left(y\right)-iS\left(y\right)=\int_{0}^{y}e^{-i\frac{\pi}{2}t^{2}}dt=y\int_{0}^{1}e^{-i\frac{\pi}{2}y^{2}s^{2}}ds=yF\left(0,i\frac{\pi}{2}y^{2}\right).\label{eq:Fresnel integrals}
\end{equation}
For computing $F\left(n_{\max},z\right)$, instead of tabulating this
function as it is done for real argument in e.g. \citep{MATSUO:1972,MAMEDO:2004,PEE-KNI:2020,WEI-OCH:2015},
we use a nonlinear approximation of the integrand in \eqref{eq:boys function}
(see \citep{BEY-MON:2016}) leading an approximation of the Boys function
valid for the complex argument $\mathcal{R}e\left(z\right)\ge0$ with
tight error estimates. For $\mathcal{R}e\left(z\right)<0$ we compute
$e^{z}F\left(n,z\right)$ instead of $F\left(n,z\right)$. Based on
these approximations, we develop two algorithms, for real and complex
valued arguments. We refer to \citep{WEIDEM:1994,CAR-POL:1998,MATHAR:2004}
for previously developed algorithms for the Boys function with complex
argument. The complex argument appears in a number of problems, for
example, in calculations with mixed Gaussian/plane wave bases in molecules
and scattering problems \citep{RE-MC-MC:1975,OSTLUN:1975,CO-FO-SI:1987,CO-FO-SI:1988,FUS-PUL:2002},
in the context of complex scaling calculations of excited states \citep{REINHA:1982},
and in using gauge invariant basis functions for calculating magnetic
properties\citep{WO-HI-PU:1990}.

\section{\label{sec:Approximation-of-F(0,z)}Approximation of $F\left(0,z\right)$
for complex valued argument}

We have 
\begin{equation}
F\left(0,z\right)=\int_{0}^{1}e^{-zt^{2}}dt=\frac{1}{2}\int_{0}^{1}e^{-zs}s^{-1/2}ds\label{eq: zero boys function}
\end{equation}
and use the integral
\begin{equation}
s^{-1/2}=\frac{2}{\sqrt{\pi}}\int_{0}^{\infty}e^{-st^{2}}dt\label{eq:integral}
\end{equation}
to obtain 
\begin{equation}
F\left(0,z\right)=\frac{1}{\sqrt{\pi}}\int_{0}^{\infty}\frac{1-e^{-\left(t^{2}+z\right)}}{t^{2}+z}dt=\frac{1}{\sqrt{\pi}}\int_{0}^{\infty}q\left(t^{2}+z\right)dt,\label{eq:final integral0}
\end{equation}
where
\[
q\left(\xi\right)=\left(1-e^{-\xi}\right)/\xi,\,\,\,\,\xi\in\mathbb{C}
\]
is an analytic function. An algorithm for computing $F\left(0,z\right)$
is essentially a quadrature for the integral in \eqref{eq:final integral0}.
Note that if, instead, we were to use a quadrature to compute $F\left(0,z\right)$
via integrals in \eqref{eq: zero boys function} then, for each $z$,
we would need to evaluate as many exponentials as the number of quadrature
terms. Importantly, when using \eqref{eq:final integral0}, we need
to evaluate $e^{-z}$ only once and then use the result as a factor.

\subsection{The case $\mathcal{R}e\left(z\right)\ge0$.}

We split the integral \eqref{eq:final integral0} into three terms
\begin{equation}
F\left(0,z\right)=\frac{1}{\sqrt{\pi}}\int_{0}^{\infty}\frac{1}{t^{2}+z}dt-\frac{e^{-z}}{\sqrt{\pi}}\int_{0}^{t_{\max}}\frac{e^{-t^{2}}}{t^{2}+z}dt-\frac{e^{-z}}{\sqrt{\pi}}\int_{t_{\max}}^{\infty}\frac{e^{-t^{2}}}{t^{2}+z}dt,\label{eq:split pos real part}
\end{equation}
and observe that the last term in \eqref{eq:split pos real part}
(without the factor $e^{-z}$) is estimated as
\begin{eqnarray}
\left|\frac{1}{\sqrt{\pi}}\int_{t_{\max}}^{\infty}\frac{e^{-t^{2}}}{t^{2}+z}dt\right| & \le & \frac{1}{\sqrt{\pi}}\int_{t_{\max}}^{\infty}\frac{e^{-t^{2}}}{\left|t^{2}+z\right|}dt\le\nonumber \\
 & \le & \frac{1}{\sqrt{\pi}}\int_{t_{\max}}^{\infty}\frac{e^{-t^{2}}}{t^{2}}dt\nonumber \\
 & = & \frac{1}{\sqrt{\pi}}\left(\frac{e^{-t_{\max}^{2}}}{t_{\max}}-\sqrt{\pi}\mbox{Erfc}\left(t_{\max}\right)\right)=\epsilon_{t_{\max}}.\label{eq:error est}
\end{eqnarray}
Selecting $t_{\max}=e^{7/4}$ to obtain $\epsilon_{t_{\max}}\approx5.9\cdot10^{-18}$.
For the first term in \eqref{eq:split pos real part} we have 
\[
\frac{1}{\sqrt{\pi}}\int_{0}^{\infty}\frac{1}{t^{2}+z}dt=\frac{1}{2}\sqrt{\frac{\pi}{z}}.
\]
For $\left|z\right|\ge r_{0}=0.35$ we use quadrature (see Appendix
for details) to approximate the second term in \eqref{eq:split pos real part}
as
\begin{equation}
\left|\frac{1}{\sqrt{\pi}}\int_{0}^{t_{\max}}\frac{e^{-t^{2}}}{t^{2}+z}dt-\sum_{m=1}^{M}\frac{w_{m}e^{-\eta_{m}}}{\eta_{m}+z}\right|\le\epsilon,\label{EQ:QUAD}
\end{equation}
where $M=22$ and nodes and weights are given in Table~\ref{tab:exp =000026 wei}.
We note that it is possible to use the standard Gauss-Legendre quadrature
on the interval $\left[0,t_{\max}\right]$ but the number of terms,
$M$, will be larger. As a result we obtain approximation

\begin{equation}
\left|F\left(0,z\right)-\left(\frac{1}{2}\sqrt{\frac{\pi}{z}}-\frac{1}{2\sqrt{\pi}}e^{-z}\sum_{m=1}^{22}\frac{w_{m}e^{-\eta_{m}}}{\eta_{m}+z}\right)\right|\le2\epsilon+\epsilon_{t_{\max}},\,\,\,\left|z\right|\ge r_{0}.\label{eq:approximation formula}
\end{equation}

\subsection{The case $\mathcal{R}e\left(z\right)<0$}

In this case, we compute $e^{z}F\left(0,z\right)$ rather than $F\left(0,z\right)$.
Since the denominator in \eqref{eq:final integral0} can be zero,
we cannot separate terms in $q\left(t^{2}+z\right)$ as in \eqref{eq:split pos real part}.
Instead, we split the integral \eqref{eq:final integral0} into two
terms and obtain
\begin{equation}
e^{z}F\left(0,z\right)=\frac{e^{z}}{\sqrt{\pi}}\int_{0}^{t_{\max}}\frac{1-e^{-\left(t^{2}+z\right)}}{t^{2}+z}dt+\frac{1}{\sqrt{\pi}}\int_{t_{\max}}^{\infty}\frac{e^{z}-e^{-t^{2}}}{t^{2}+z}dt.\label{eq:two terms}
\end{equation}
The first term in \eqref{eq:two terms} is approximated by using the
Gauss-Legendre quadrature on the interval $\left[0,t_{\max}\right]$.
The function $q$ is analytic and, therefore, there is no singularity
at $t^{2}=-z$. Since we can compute derivatives of $q$, the error
introduced by this quadrature can be estimated using results in \citep[Section 5.2]{KA-MO-NA:1989}.
For example, we obtain 
\[
\left|\frac{1}{\sqrt{\pi}}\int_{0}^{t_{\max}}\frac{e^{z}-e^{-t^{2}}}{t^{2}+z}dt-\frac{1}{\sqrt{\pi}}\sum_{m=1}^{M^{g}}w_{m}^{g}\frac{e^{z}\left(1-e^{-\left(t_{m}^{2}+z\right)}\right)}{t_{m}^{2}+z}\right|\le\epsilon^{g}
\]
with $M^{g}=16$ and $\epsilon^{g}\approx10^{-14}$ where $t_{m,}w_{m}^{g}$
are the standard Gauss-Legendre nodes and weights on the interval
$\left[0,t_{\max}\right]$. In the second term in \eqref{eq:two terms},
we drop $e^{-t^{2}}$(since its contribution is less than $e^{-t_{max}^{2}}\approx4.2\cdot10^{-15}$)
and obtain
\begin{equation}
\frac{1}{\sqrt{\pi}}\int_{t_{\max}}^{\infty}\frac{e^{z}}{t^{2}+z}dt=\frac{e^{z}}{\sqrt{\pi z}}\mbox{Arctan}\left(\frac{\sqrt{z}}{t_{\max}}\right).\label{eq:approx with arctan}
\end{equation}
While we obtain an explicit expression, computing arctangent of a
complex argument is relatively expensive. For a complex argument,
we evaluate arctangent using 
\[
\mbox{Arctan}\left(z\right)=\frac{1}{2}i\log\frac{1-iz}{1+iz}.
\]
As a result of dropping $e^{-t^{2}}$ in the second term of \eqref{eq:two terms},
our approximation in \eqref{eq:approx with arctan} has a singularity
at $z=-t_{\max}^{2}$. In order to avoid using \eqref{eq:approx with arctan}
in the vicinity of singularity, we use two different parameters, $t_{\max}$
and $t_{\max,1}$ and switch to the version with $t_{\max,1}$ if
$\left|z+t_{\max}^{2}\right|\le1/2$, where $t_{\max,1}=\sqrt{t_{\max}^{2}+1}$.

We note that it is possible to increase the number of terms in the
quadrature in order to avoid evaluating arctangent. This might be
of interest on a parallel (GPU or multi-core) computer since computation
of quadrature terms is trivially parallel. As a result, we obtain
approximation
\begin{equation}
\left|e^{z}F\left(0,z\right)-\frac{e^{z}}{\sqrt{\pi z}}\mbox{Arctan}\left(\frac{\sqrt{z}}{t_{\max}}\right)-\frac{1}{\sqrt{\pi}}\sum_{m=1}^{M^{g}}w_{m}^{g}\frac{e^{z}\left(1-e^{-\left(t_{m}^{2}+z\right)}\right)}{t_{m}^{2}+z}\right|\le\tilde{\epsilon},\,\,\,\left|z+t_{\max}^{2}\right|>1/2,\label{eq:approx for negative real part}
\end{equation}
where $\tilde{\epsilon}\approx10^{-14}$. For $\left|z\right|>t_{\max}^{2}$,
we have a converging series for the second integral in \eqref{eq:split pos real part}
as follows:
\begin{eqnarray}
\int_{0}^{t_{\max}}\frac{e^{-t^{2}}}{t^{2}+z}dt & = & \frac{1}{z}\int_{0}^{t_{\max}}\frac{e^{-t^{2}}}{t^{2}/z+1}dt\nonumber \\
 & = & \frac{1}{z}\sum_{j=0}^{\infty}\left(-1\right)^{j}z^{-j}\int_{0}^{t_{\max}}e^{-t^{2}}t^{2j}dt\nonumber \\
 & = & \frac{1}{2z}\sum_{j=0}^{\infty}\left(-1\right)^{j}z^{-j}\left(\Gamma\left(j+1/2\right)-\Gamma\left(j+1/2,t_{\max}^{2}\right)\right)\label{eq:series}\\
 & = & \frac{1}{z}\sum_{j=0}^{\infty}\left(-1\right)^{j}z^{-j}t_{\max}^{2j+1}F\left(j,t_{\max}^{2}\right)\nonumber 
\end{eqnarray}
so that we can use 
\begin{equation}
\left|F\left(0,z\right)-\left[\frac{\sqrt{\pi}}{2\sqrt{z}}-\frac{e^{-z}}{2\sqrt{\pi}z}\sum_{j=0}^{J}\left(-1\right)^{j}z^{-j}\left(\Gamma\left(j+1/2\right)-\Gamma\left(j+1/2,t_{\max}^{2}\right)\right)\right]\right|\le\epsilon_{t_{\max}},\,\,\,\,\left|z\right|>t_{\max}^{2},\label{eq:large arg F_0}
\end{equation}
instead of \eqref{eq:approximation formula} and 
\begin{equation}
\left|e^{z}F\left(0,z\right)-\left[\frac{e^{z}\sqrt{\pi}}{2\sqrt{z}}-\frac{1}{2\sqrt{\pi}z}\sum_{j=0}^{J}\left(-1\right)^{j}z^{-j}\left(\Gamma\left(j+1/2\right)-\Gamma\left(j+1/2,t_{\max}^{2}\right)\right)\right]\right|\le\epsilon_{t_{\max}},\,\,\,\,\left|z\right|>t_{\max}^{2},\label{eq:large arg F_0-1}
\end{equation}
instead of \eqref{eq:approx for negative real part}. Since the parameter
$t_{\max}$ is fixed, the coefficients of the series are computed
offline.

Note that the series in \eqref{eq:series}-\eqref{eq:large arg F_0}
is related to the asymptotic expansion of $F\left(0,z\right)$ (see
e.g. \citep{CAR-POL:1998}),
\begin{eqnarray}
F\left(0,z\right) & \sim & \frac{1}{2}\frac{\sqrt{\pi}}{\sqrt{z}}-\frac{\sqrt{\pi}}{2}\frac{e^{-z}}{z}\sum_{j=0}^{J}\frac{z^{-j}}{\Gamma\left(\frac{1}{2}-j\right)}\label{eq:asymptotics F_0}\\
 & = & \frac{1}{2}\frac{\sqrt{\pi}}{\sqrt{z}}-\frac{e^{-z}}{2\sqrt{\pi}z}\sum_{j=0}^{J}\left(-1\right)^{j}z^{-j}\Gamma\left(j+1/2\right).\nonumber 
\end{eqnarray}
We use \eqref{eq:large arg F_0} and \eqref{eq:large arg F_0-1} for
$\left|z\right|\ge100$ so that it is sufficient to keep only seven
terms yielding an error of less than $10^{-13}$.

For $\left|z\right|\le r_{0}$ we use the Taylor expansion of \eqref{eq: zero boys function},
\begin{equation}
F\left(0,z\right)=\sum_{j=0}^{\infty}\frac{\left(-1\right)^{j}z^{j}}{j!\left(2j+1\right)}.\label{eq:Taylor for F_0}
\end{equation}
and we need $10$ terms to maintain an accuracy of about $13$ digits.
While selecting parameters as above leads to algorithms with a reasonable
speed, we did not optimize these choices as they may depend on several
factors, e.g. computer architecture.

Since the Boys function $F\left(0,z\right)$ is related to the error
function (and can be used to compute it), we compared the speed of
our algorithm with that of the well-known algorithm by Gautschi \citep{GAUTSC:1970}
for computing the error function with complex argument using a rational
approximation of the closely related Faddeeva function. The speed
of that algorithm was measured in comparison with the speed of computing
$\exp\left(z\right)$. In \citep{GAUTSC:1970} it is stated that with
an accuracy of $\sim10$ digits, the code is $7-15$  times slower
than the speed of computing $\exp\left(z\right)$. Using the same
comparison for our algorithm, this ratio is  $\sim12$ for an accuracy
of about $13$ digits. Our algorithm is implemented using Fortran
90 compiled by Intel's ifort with compiler flags -O3 -ipo -static
and running on a laptop with $\approx2.3$ GHz chipset. We timed our
code by performing $10^{6}$ evaluations yielding $\approx0.92\cdot10^{-7}$
seconds per evaluation in comparison with $\approx0.79\cdot10^{-8}$
seconds per evaluation for $\exp\left(z\right)$ with a complex argument.

While algorithms for computing the Fresnel integrals appear to be
somewhat faster than using the Boys function in \eqref{eq:Fresnel integrals}
(see e.g. \citep{BOERSM:1960}), we note that the generalized Fresnel
integrals, e.g. $\int_{0}^{x}e^{it^{n}}dt$, $n\ge2$, can be evaluated
using our approach and plan to consider algorithms for these oscillatory
special functions elsewhere.

\begin{table}[H]
\begin{centering}
\begin{tabular}{|c|l|l|c|l|l|}
\hline 
$m$ & $\eta_{m}$ & $w_{m}\cdot e^{-\eta_{m}}$ & $m$ & $\eta_{m}$ & $w_{m}\cdot e^{-\eta_{m}}$\tabularnewline
\hline 
$1$ & $0.14778782637969565E-02$ & $0.86643102720141654E-01$ & $12$ & $0.12539502287919293E+01$ & $0.57444804221430223E-01$\tabularnewline
\hline 
$2$ & $0.13317276413725817E-01$ & $0.85772060843439468E-01$ & $13$ & $0.17244634233573395E+01$ & $0.42081994346945442E-01$\tabularnewline
\hline 
$3$ & $0.37063591452052541E-01$ & $0.83935043682917876E-01$ & $14$ & $0.23715248262781863E+01$ & $0.25838539448223272E-01$\tabularnewline
\hline 
$4$ & $0.72752512422882762E-01$ & $0.80966197041322921E-01$ & $15$ & $0.32613796996078355E+01$ & $0.12445024157255560E-01$\tabularnewline
\hline 
$5$ & $0.12023694122878568E+00$ & $0.76908954849297856E-01$ & $16$ & $0.44851301690595911E+01$ & $0.42925415925998368E-02$\tabularnewline
\hline 
$6$ & $0.17957429395893773E+00$ & $0.73155207871182168E-01$ & $17$ & $0.61680621351224838E+01$ & $0.93543429877359686E-03$\tabularnewline
\hline 
$7$ & $0.25353404698408727E+00$ & $0.72695003516315720E-01$ & $18$ & $0.84824718723178698E+01$ & $0.10840885466502505E-03$\tabularnewline
\hline 
$8$ & $0.35038865278072195E+00$ & $0.75284255608930400E-01$ & $19$ & $0.11665305486296793E+02$ & $0.52718679667616736E-05$\tabularnewline
\hline 
$9$ & $0.48210957593127668E+00$ & $0.77094395364519633E-01$ & $20$ & $0.16042417132288328E+02$ & $0.77659740397504190E-07$\tabularnewline
\hline 
$10$ & $0.66302899315837416E+00$ & $0.75425062567753040E-01$ & $21$ & $0.22061929518147089E+02$ & $0.22138172422680093E-09$\tabularnewline
\hline 
$11$ & $0.91181473685659087E+00$ & $0.68968619265031533E-01$ & $22$ & $0.30340112094708307E+02$ & $0.65941617600377069E-13$\tabularnewline
\hline 
\end{tabular}
\par\end{centering}
\caption{\label{tab:exp =000026 wei}The poles and weights in \eqref{eq:approximation formula}.}
\end{table}

\section{Approximation of $F\left(n_{\max},z\right)$ for real and complex
arguments}

The function 
\begin{equation}
g_{n}\left(s\right)=\left(1-s\right){}^{n-1/2}\label{eq:function g}
\end{equation}
decays monotonically on $\left[0,1\right]$ and we use Algorithm~1
in \citep{BE-MO-SA:2018} to construct its near optimal approximation
via exponentials. We refer to \citep{BE-MO-SA:2018} and references
therein for the details of the algorithm that we use to obtain the
necessary parameters (for an example, see Table~\ref{tab:Weights-and-exponents}).

We obtain approximation
\begin{equation}
\left|g_{n}\left(s\right)-\sum_{m=1}^{M}w_{m}e^{\eta_{m}s}\right|\le\epsilon\,\,\,\,\,\mbox{for}\,\,\,s\in\left[0,1\right],\label{eq:key approximation}
\end{equation}
where $w_{m},\eta_{m}\in\mathbb{C}$. We note that $n$ should be
sufficiently large (e.g. $n\ge7)$ to avoid the impact on the approximation
of the singularity of the $n$-th derivative of $g_{n}$. Its numerical
effect makes the accuracy of the current double precision implementation
of Algorithm~1 in \citep{BE-MO-SA:2018} insufficient to reliably
produce approximation \eqref{eq:key approximation} for $1\le n\le6$.

Substituting the approximation of $g_{n}\left(1-s\right)$ into the
integral \eqref{eq:boys function}, we arrive at 
\[
F\left(n,z\right)-\frac{1}{2}\sum_{m=1}^{M}w_{m}\int_{0}^{1}e^{-zs}e^{\eta_{m}\left(1-s\right)}ds=\frac{1}{2}\int_{0}^{1}e^{-zs}\left(s^{n-1/2}-\sum_{m=1}^{M}w_{m}e^{\eta_{m}\left(1-s\right)}\right)ds
\]
and estimate
\[
\left|F\left(n,z\right)-\frac{1}{2}\sum_{m=1}^{M}w_{m}\int_{0}^{1}e^{-zs}e^{\eta_{m}\left(1-s\right)}ds\right|\le\frac{\epsilon}{2}\int_{0}^{1}e^{-\mathcal{R}e\left(z\right)s}ds=\frac{\epsilon}{2}\frac{1-e^{-\mathcal{R}e\left(z\right)}}{\mathcal{R}e\left(z\right)}.
\]
Since
\[
\frac{1}{2}\sum_{m=1}^{M}w_{m}\int_{0}^{1}e^{-zs}e^{\eta_{m}\left(1-s\right)}ds=\frac{1}{2}\sum_{m=1}^{M}w_{m}e^{\eta_{m}}\frac{1-e^{-\left(z+\eta_{m}\right)}}{z+\eta_{m}},
\]
we have
\begin{equation}
\left|F\left(n,z\right)-\frac{1}{2}\sum_{m=1}^{M}w_{m}e^{\eta_{m}}\frac{1-e^{-\left(z+\eta_{m}\right)}}{z+\eta_{m}}\right|\le\frac{\epsilon}{2}\frac{1-e^{-\mathcal{R}e\left(z\right)}}{\mathcal{R}e\left(z\right)}\le\frac{\epsilon}{2}.\label{eq:estimate}
\end{equation}
Indeed, denoting the factor on the right hand side of \eqref{eq:estimate},
$q\left(z\right)=\left(1-e^{-z}\right)/z$, we have 
\[
q\left(2z\right)=\frac{1}{2}\left(e^{-z}q\left(z\right)+q(z)\right)
\]
and, therefore, for $\mathcal{R}e\left(z\right)\ge0$
\[
\left|q\left(2z\right)\right|\le\left|q\left(z\right)\right|.
\]
This implies that $\left|q\left(z\right)\right|$ reaches it maximum
at $z=0$ , where $q\left(0\right)=1$.

If $\mathcal{R}e\left(z\right)<0$ we compute $e^{z}F\left(0,z\right)$
instead of $F\left(0,z\right)$,
\[
e^{z}F\left(n,z\right)=\frac{1}{2}\int_{0}^{1}e^{z(1-s)}s^{n-1/2}ds=\frac{1}{2}\int_{0}^{1}e^{zs}\left(1-s\right)^{n-1/2}ds.
\]
Using \eqref{eq:key approximation} we obtain
\[
e^{z}F\left(n,z\right)-\frac{1}{2}\sum_{m=1}^{M}w_{m}\int_{0}^{1}e^{zs}e^{\eta_{m}s}ds=\frac{1}{2}\int_{0}^{1}e^{zs}\left[g_{n}\left(s\right)-\sum_{m=1}^{M}w_{m}e^{\eta_{m}s}\right]ds
\]
and the estimate
\[
\left|e^{z}F\left(n,z\right)-\frac{1}{2}\sum_{m=1}^{M}w_{m}\frac{e^{z+\eta_{m}}-1}{z+\eta_{m}}\right|\le\frac{\epsilon}{2}\frac{e^{\mathcal{R}e\left(z\right)}-1}{\mathcal{R}e\left(z\right)}\le\frac{\epsilon}{2}.
\]
For computing values of $e^{z}F\left(n,z\right)$ for $0\le n\le n_{\max}$
for $\mathcal{R}e\left(z\right)<0$, we use recursions
\begin{equation}
e^{z}F\left(n,z\right)=\frac{n-1/2}{z}e^{z}F\left(n-1,z\right)-\frac{1}{2z}\label{eq:forward modified}
\end{equation}
instead of \eqref{eq:forward recursion} and 
\begin{equation}
e^{z}F\left(n-1,z\right)=\frac{2x}{2n-1}e^{z}F\left(n,z\right)+\frac{1}{2n-1}\label{eq:backward modified}
\end{equation}
instead of \eqref{eq:backward recursion}.

\section{Implementation details}

\begin{figure}[H]
\begin{centering}
\includegraphics[scale=0.7]{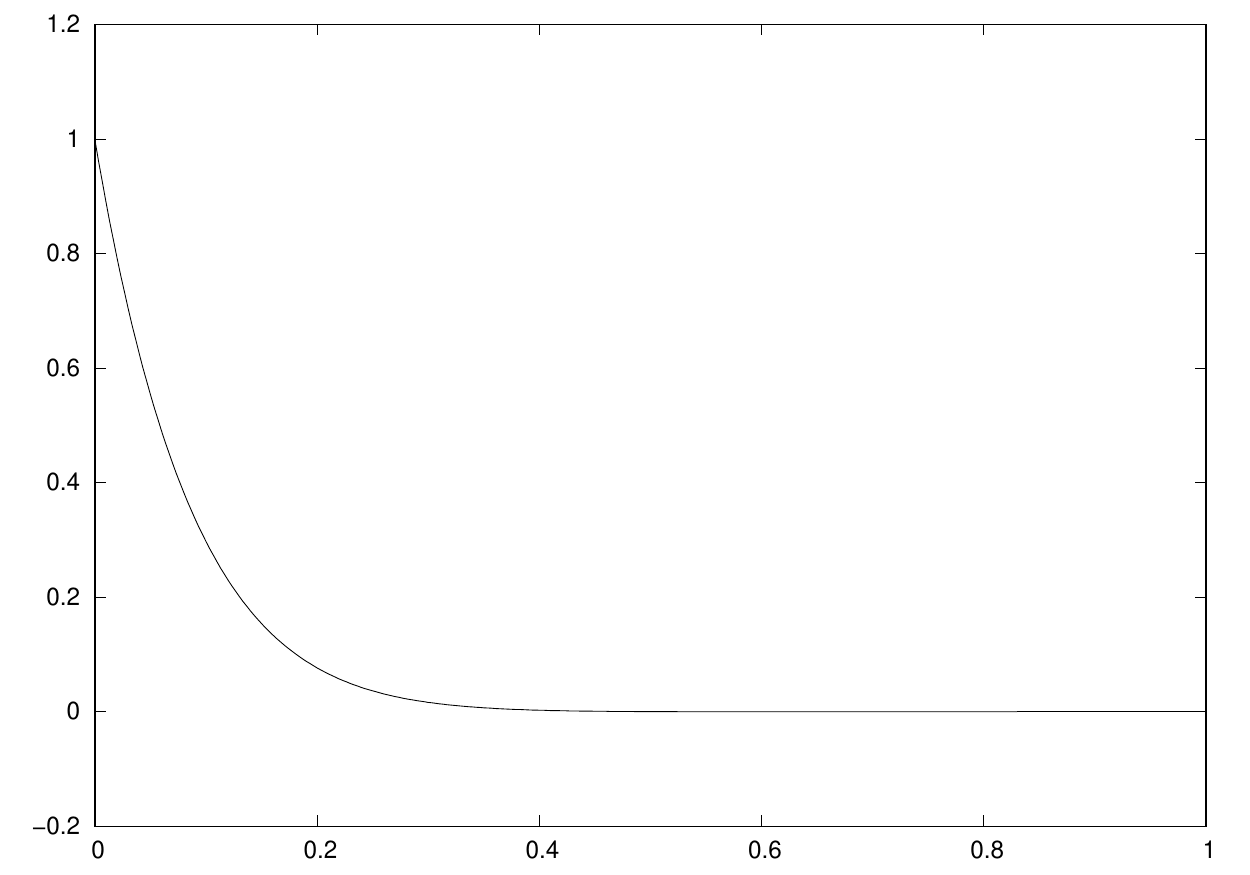}\includegraphics[scale=0.7]{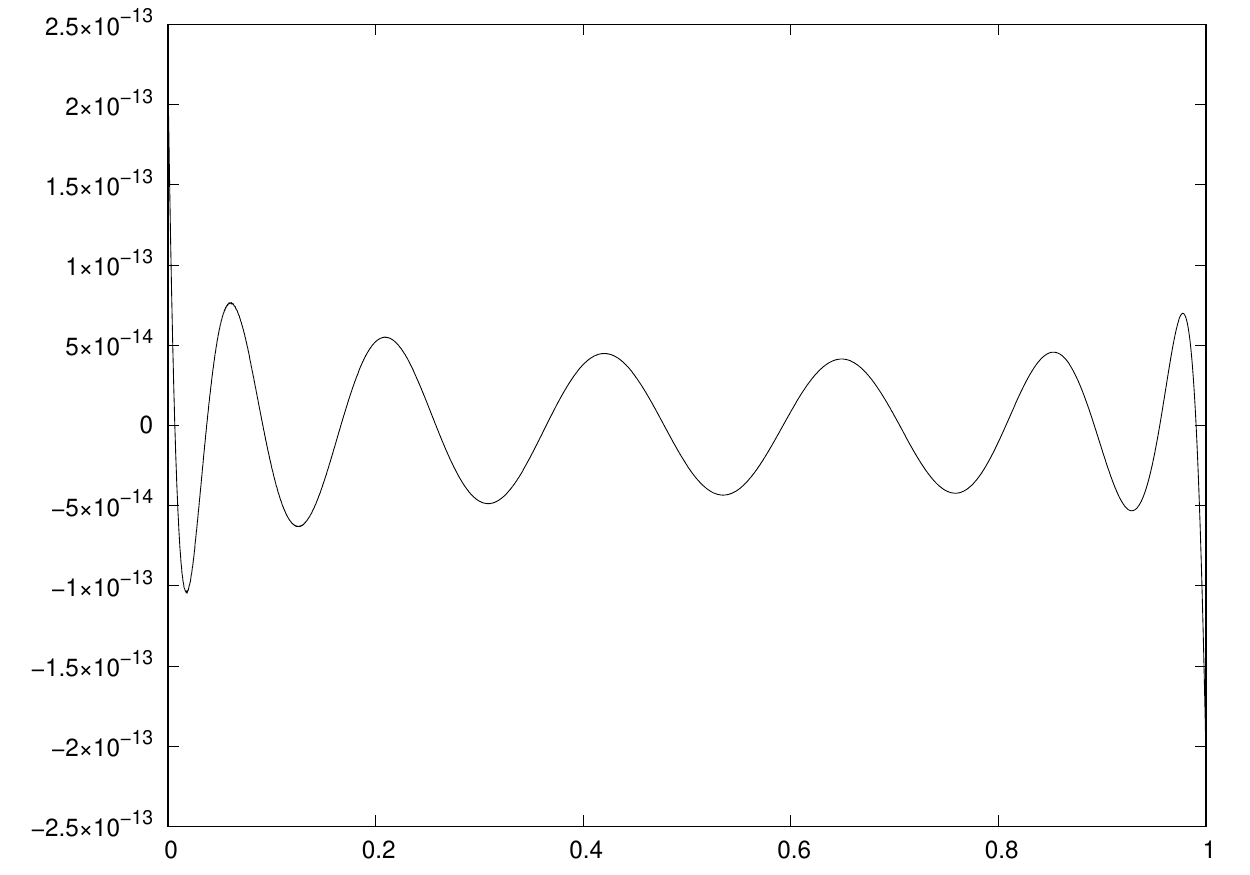}
\par\end{centering}
\caption{\label{fig:Function--in}The function $g_{12}\left(s\right)$ in \eqref{eq:function g}
and the error of its near optimal approximation via exponentials in
\eqref{eq:key approximation} with parameters described in Table~\ref{tab:Weights-and-exponents}.}
\end{figure}

The speed of computation of values of $F\left(n_{\max},z\right)$
for $n_{\max}\ge7$ depends on the number of terms $M$ in approximation
\eqref{eq:key approximation}. We demonstrate the results of approximating
$F\left(12,z\right)$ and display function $g_{12}\left(s\right)$
in Figure~\ref{fig:Function--in}. Using only 13 terms (see Table~\ref{tab:Weights-and-exponents}),
we achieve accuracy for $F\left(12,z\right)$ $\epsilon\approx2\cdot10^{-14}$
(e.g. accuracy of evaluation of $F\left(12,0\right)$ is $2.08\cdot10^{-14}$).

In implementing this approximation, we need to isolate cases where
$z$ is close to $-\eta_{m}$ by using the Taylor expansion for $\frac{1-e^{-\left(z+\eta_{m}\right)}}{z+\eta_{m}}$.
Since most of $\eta_{m}$ have imaginary part, it is a minimal effort
if $z$ is real since $\eta_{m}$ is real in only three terms in our
example in Table~\ref{tab:Weights-and-exponents}. In addition, for
the real argument $z$, we need to use only five terms with complex
valued parameters as they come in complex conjugate pairs.

We implemented these algorithms using Fortran 90 on a laptop described
in Section~\ref{sec:Approximation-of-F(0,z)}. Computing the Boys
functions $F\left(n,z\right)$ for $n=0,\dots12$ for the real argument
takes $\approx0.34\cdot10^{-7}$ seconds. The subroutine for the complex
valued argument is slower and takes $\approx0.21\cdot10^{-6}$ seconds.

\begin{table}[h]
\begin{centering}
\begin{tabular}{|c|l|l|}
\hline 
$m$ & $\eta_{m}$ & $w_{m}$\tabularnewline
\hline 
\hline 
$1$ & $0.70719431320570010\cdot10^{1}+0.16487291250752115\cdot10^{2}i$ & $0.36443632402898501\cdot10^{-10}+0.26411751072107504\cdot10^{-10}i$\tabularnewline
\hline 
$2$ & $0.70719431320570010\cdot10^{1}-0.16487291250752115\cdot10^{2}i$ & $0.36443632402898501\cdot10^{-10}-0.26411751072107504\cdot10^{-10}i$\tabularnewline
\hline 
$3$ & $-0.57143271715191635+0.13278579453233633\cdot10^{2}i$ & $0.18185250346753633\cdot10^{-6}-0.21860458971399352\cdot10^{-5}i$\tabularnewline
\hline 
$4$ & $-0.57143271715191635-0.13278579453233633\cdot10^{2}i$ & $0.18185250346753633\cdot10^{-6}+0.21860458971399352\cdot10^{-5}i$\tabularnewline
\hline 
$5$ & $-0.47193021330392506\cdot10^{1}+0.99835257112371032\cdot10^{1}i$ & $-0.99489169272055748\cdot10^{-3}-0.23049079105203073\cdot10^{-3}i$\tabularnewline
\hline 
$6$ & $-0.47193021330392506\cdot10^{1}-0.99835257112371032\cdot10^{1}i$ & $-0.99489169272055748\cdot10^{-3}+0.23049079105203073\cdot10^{-3}i$\tabularnewline
\hline 
$7$ & $-0.71704662772895089\cdot10^{1}+0.66712360839820768\cdot10^{1}i$ & $-0.25625216985879006\cdot10^{-1}+0.35818335274876982\cdot10^{-1}i$\tabularnewline
\hline 
$8$ & $-0.71704662772895089\cdot10^{1}-0.66712360839820768\cdot10^{1}i$ & $-0.25625216985879006\cdot10^{-1}-0.35818335274876982\cdot10^{-1}i$\tabularnewline
\hline 
$9$ & $-0.84899747054724699\cdot10^{1}+0.33434804168467491\cdot10^{1}i$ & $0.16506801544880723+0.32273964471776045i$\tabularnewline
\hline 
$10$ & $-0.84899747054724699\cdot10^{1}-0.33434804168467491\cdot10^{1}i$ & $0.16506801544880723-0.32273964471776045i$\tabularnewline
\hline 
$11$ & $0.36564414363150973\cdot10^{2}$ & $-0.20104641661565164\cdot10^{-25}$\tabularnewline
\hline 
$12$ & $-0.32424239255921954\cdot10^{1}$ & $-0.39563536955042078\cdot10^{-3}$\tabularnewline
\hline 
$13$ & $-0.89066047733100753\cdot10^{1}$ & $0.72349945805085292$\tabularnewline
\hline 
\end{tabular}
\par\end{centering}
\caption{\label{tab:Weights-and-exponents}Weights and exponents of the approximation
of $g_{12}\left(s\right)$ on $\left[0,1\right]$ in \eqref{eq:key approximation}.
With these parameters the absolute error in \eqref{eq:estimate} is
$\epsilon\approx2.5\cdot10^{-13}$ .}
\end{table}

\section{Conclusion}

Since their introduction in \citep{BOYS:1950}, the Boys functions
with real argument have widely been  used for computing Gaussian integrals.
When using mixed Gaussian/exponential bases, one needs to evaluate
the Boys functions with complex argument. Such mixed bases are appropriate
for scattering problems and for bound state problems where using only
plane waves becomes too expensive near singularities. Consequently,
mixed Gaussian/exponential bases provide a greater flexibility in
formulation and solving problems of quantum chemistry and we present
our results in part to facilitate their use.

While for real argument the Boys functions can be easily tabulated
in regions where their asymptotic is not accurate, it is more difficult
to apply such straightforward implementation for a complex argument.
A careful reading of references \citep{WEIDEM:1994,CAR-POL:1998,MATHAR:2004}
reveals shortcomings of existing approaches (relying mostly on expansions)
to computing the Boys functions of complex argument (see e.g. conclusion
in \citep{MATHAR:2004}). For our approach a better comparison is
offered by Gautschi's algorithm \citep{GAUTSC:1970} for the error
function of complex argument since it is related to $F\left(0,z\right)$
as in \eqref{eq:Boys and Error fnc}, see Section~\ref{sec:Approximation-of-F(0,z)}.
Our approach of approximating a part of the integrand so that the
resulting integral can be evaluated explicitly, is simpler and yields
tight accuracy estimates. Note that the part of the integrand we are
approximating is real while the Boys functions we are computing are
complex-valued. As a side remark we note that the Boys function $F\left(0,z\right)$
remains bounded for complex argument with $\mathcal{R}e\left(z\right)\ge0$
(and $e^{z}F\left(0,z\right)$ for $\mathcal{R}e\left(z\right)<0$)
and, for this reason, provides a good alternative approach for computing
the error function of complex argument.

We avoid the direct timing comparisons with existing algorithms since
such comparisons are generally misleading. Given different hardware
(single core, multi-core, GPU, etc), different compilers and compiler
flags, and different implementations, it is hard to compare algorithms
by simply running them. Instead one can look at algorithmic possibilities
an approach offers. Our code is compact and it is easy to simply count
the total number of operations. We note that computation of each term
in the sums \eqref{eq:approximation formula} and \eqref{eq:estimate}
is trivially parallel and only recursions in \eqref{eq:forward recursion}
and \eqref{eq:backward recursion} require a sequential implementation
(with just three multiplications and one addition per function). Thus
timing of our algorithms implemented on a multi-core or GPU computer
will be much faster than the quoted timings of our implementation
on a single CPU.

\section{Supplementary Material}

The supplementary material for this paper consists of $5$ Fortran
90 subroutines implementing, as an example, algorithms for computing
the Boys function with indices $n=0,\dots12$. The subroutine dboysfun12.f90
evaluates the Boys functions $F\left(n,z\right)$ for real non-negative
argument $z$. The subroutines zboysfun12.f90 and zboysfun00.f90 evaluate
the Boys functions $F\left(n,z\right)$ for complex argument $z$
with non-negative real part. Finally, the subroutines zboysfun00nrp.f90,
zboysfun12nrp.f90 evaluate the functions $e^{z}F\left(n,z\right)$
for complex argument $z$ with negative real part.

\section{Acknowledgements}

SS was supported by the NSF under Grant CHE-1800584. SS was also partly
supported through the Sloan research fellowship.

\section{Data availability}

The data that supports the findings of this study are available within
the article and its supplementary material.

\section{Appendix: Construction of quadrature in \eqref{EQ:QUAD}}

Changing variables in \eqref{eq:integral} $t=e^{\tau/2}$, we rewrite
it as 
\begin{equation}
s^{-1/2}=\frac{1}{\sqrt{\pi}}\int_{-\infty}^{\infty}e^{-se^{\tau}+\tau/2}d\tau,\,\,\,\,0\le s\le1,\label{eq:integral-2}
\end{equation}
 and discretize it (following \citep{BEY-MON:2010}) yielding

\begin{equation}
\left|s^{-1/2}-\sum_{m=1}^{M}w_{m}e^{-\eta_{m}s}\right|\le\epsilon s^{-1/2},\,\,\,\delta\le s\le1,\label{eq:approx s^(-1/2)}
\end{equation}
where $\eta_{m},w_{m}>0$ are arranged in an ascending order, and
we estimate that
\begin{eqnarray*}
\left|F\left(0,z\right)-\frac{1}{2}\int_{0}^{1}e^{-zs}\left(\sum_{m=1}^{M}w_{m}e^{-\eta_{m}s}\right)ds\right| & \le & \frac{1}{2}\int_{0}^{1}e^{-\mathcal{R}e\left(z\right)s}\left|s^{-1/2}-\sum_{m=1}^{M}w_{m}e^{-\eta_{m}s}\right|ds\\
 & \le & \frac{\epsilon}{2}\int_{0}^{1}e^{-z\mathcal{R}e\left(z\right)}s^{-1/2}ds\\
 & = & \epsilon F\left(0,\mathcal{R}e\left(z\right)\right)\le\epsilon.
\end{eqnarray*}
Using $\delta=\epsilon=10^{-13}$ in \eqref{eq:approx s^(-1/2)} results
in approximation with $M=210$. We also need this approximation to
satisfy
\begin{equation}
\left|\frac{1}{2\sqrt{\pi}}\int_{-\infty}^{\infty}\frac{1}{e^{\tau}+z}e^{\tau/2}d\tau-\sum_{m=1}^{M}\frac{w_{m}}{\eta_{m}+z}\right|=\left|\frac{1}{2}\sqrt{\frac{\pi}{z}}-\sum_{m=1}^{M}\frac{w_{m}}{\eta_{m}+z}\right|\le\epsilon,\,\,\,\,\left|z\right|\ge r_{0},\label{eq:prelim0}
\end{equation}
in order to obtain
\begin{equation}
\left|F\left(0,z\right)-\left(\frac{1}{2}\sqrt{\frac{\pi}{z}}-\frac{e^{-z}}{2\sqrt{\pi}}\sum_{m=1}^{M}\frac{w_{m}e^{-\eta_{m}}}{\eta_{m}+z}\right)\right|\le2\epsilon.\label{eq:prelim expression}
\end{equation}
The exponents and the weights in \eqref{eq:approx s^(-1/2)} grow
as $\eta_{m}\approx e^{\tau_{m}}$ and $w_{m}\approx e^{\tau_{m}/2}$
(see \citep{BEY-MON:2010}), so that in \eqref{eq:prelim expression}
it is sufficient to use a subset of terms with $\eta_{m}\le e^{\tau_{max}}$.
Selecting $\tau_{\max}=7/2$ so that $t_{\max}=e^{\tau_{\max}/2}$
in \eqref{eq:error est}, the error $\epsilon_{t_{\max}}\approx5.9\cdot10^{-18}$.
Consequently, we only need the $22$ terms displayed in Table~\ref{tab:exp =000026 wei}
and obtain approximation of \eqref{eq:final integral0} in \eqref{eq:approximation formula}.

\bibliographystyle{unsrt}

\end{document}